\input amsppt.sty
\magnification=\magstep1
\hsize=30truecc
\vsize=22.2truecm
\baselineskip=16truept
\nologo
\pageno=1
\topmatter
\def\N{\Bbb N}
\def\Z{\Bbb Z}
\def\Q{\Bbb Q}

\def\C{\Bbb C}
\def\l{\left}
\def\r{\right}
\def\b{\bigg}
\def\bg{\bigg}
\def\({\b(}
\def\[{\b[}
\def\){\b)}
\def\]{\b]}

\def\t{\text}
\def\f{\frac}
\def\mo{\roman{mod}}
\def\em{\emptyset}
\def\se {\subseteq}

\def\sm{\setminus}

\def\bi{\binom}
\def\eq{\equiv}
\def\cs{\cdots}
\def\ls{\leqslant}
\def\gs{\geqslant}
\def\al{\alpha}

\def\Proof{\noindent{\it Proof}}
\def\Remark{\noindent{\it Remark}}
\def\endrem{\medskip}

\hbox {New version (2004-04-15), {\tt arXiv:math.NT/0403271}.}
\bigskip
\title On $m$-covers and $m$-systems\endtitle
\author Zhi-Wei Sun\endauthor
\affil Department of Mathematics, Nanjing University
\\ Nanjing 210093, The People's Republic of China\\
{\it E-mail}: {\tt zwsun$\@$nju.edu.cn}
\\Homepage: {\tt http://pweb.nju.edu.cn/zwsun}\endaffil

\abstract Let $\Cal A=\{a_s(\mo\ n_s)\}_{s=0}^k$ be a system of
residue classes. With the help of cyclotomic fields we obtain a theorem which
unifies several previously known results concerning system $\Cal A$.
In particular, we show that if every integer lies in more than
$m_0=\lfloor\sum_{s=1}^{k}1/n_s\rfloor$ members of $\Cal A$, then for any
$a=0,1,2,\ldots$ there are at least
$\bi {m_0}{\lfloor a/n_0\rfloor}$
subsets $I$ of $\{1,\ldots,k\}$ with $\sum_{s\in I}1/n_s=a/n_0$.
We also characterize when any integer lies in at most $m$ members
of $\Cal A$, where $m$ is a fixed positive integer.
\endabstract
\thanks 2000 {\it Mathematics Subject Classification}. Primary 11B25;
Secondary 05A05, 11A07, 11B75, 11D68.
\newline\indent
The author is supported by the National Science Fund for
Distinguished Young Scholars (No. 10425103) and the Key Program of
NSF (No. 10331020) in P. R. China.
\endthanks
\endtopmatter
\document
\heading {1. The main results}\endheading

  For $a\in\Z$ and $n\in\Z^+=\{1,2,3,\ldots\}$, we simply denote
the residue class $\{x\in\Z:\, x\eq a\ (\mo\ n)\}$ by $a(n)$.
For a finite system
$$A=\{a_s(n_s)\}_{s=1}^k\tag1.1$$
of residue classes,
the function $w_{A}:\Z\to \N=\{0,1,2,\ldots\}$  given by
$$ w_{A}(x)=|\{1\ls s\ls k:\,   x\in a_s(n_s)\}|\tag1.2$$
is called the {\it covering function} of $A$. Obviously $w_A(x)$
is periodic modulo the least common multiple $N$ of the moduli
$n_1,\ldots,n_k$, and it is easy to see that the average
$\sum_{x=0}^{N-1}w_A(x)/N$ equals $\sum_{s=1}^k1/n_s$. As in [S97, S99] we call
$m(A)=\min_{x\in\Z}w_A(x)$ the {\it covering multiplicity} of
system (1.1).

  Let $m$ be any positive integer.
If $w_A(x)\gs m$ for all $x\in\Z$ (i.e., $m(A)\gs m$), then $(1.1)$
is said to be an {\it $m$-cover} of
$\Z$ as in [S95, S96], and in this case $\sum_{s=1}^k1/n_s\gs m$.
The concept of cover (i.e. $1$-cover) was introduced by P. Erd\H
os [E50]. If $w_A(x)=m$ for all $x\in\Z$, then we call $(1.1)$ an
{\it exact $m$-cover} of $\Z$ as in [S96, S97] (and in this case
$\sum_{s=1}^k1/n_s=m$); such covers were first studied by \v S.
Porubsk\'y [P76]. If $w_A(x)\ls m$ for all $x\in\Z$, then we call
(1.1) an {\it $m$-system}, and in this case $\sum_{s=1}^k1/n_s\ls
m$; any $1$-system is said to be {\it disjoint}. The reader is
referred to [PS] for a survey of modern results in covering theory.
Covers of $\Z$ have many surprising applications (cf. [C71], [F02],
[S00], [S01], [S04] and [SY]);  recently the author [S03b] discovered that they
are also related to zero-sum problems and subset sums in combinatorial number theory.

 Throughout this paper, for $a,b\in\Z$ we set
  $[a,b]=\{x\in\Z:\, a\ls x\ls b\}$ and define $[a,b)$ and $(a,b]$ similarly.
 As usual, the integral part and the fractional part of a real number $\al$
 are denoted by $\lfloor\al\rfloor$ and $\{\al\}$ respectively.

  For system (1.1) we define its {\it dual system} $A^*$ by
$$A^*=\{a_s+r(n_s):\,  \ 1\ls r<n_s,\ 1\ls s\ls k\}.\tag1.3$$
As $\{a_s+r(n_s)\}_{r=0}^{n_s-1}$ is a partition of $\Z$ for any $s\in[1,k]$,
we have $w_A(x)+w_{A^*}(x)=k$ for all $x\in\Z$.
Thus $w_A(x)\ls m$ for all $x\in\Z$ if and only if
$w_{A^*}(x)\gs k-m$ for all $x\in\Z$.
This simple and new observation shows that we can study $m$-systems
via covers of $\Z$, and construct covers of $\Z$ via $m$-systems.

 By a result in [S96], if (1.1) is an $m$-cover of $\Z$
then for any $m_1,\ldots,m_k\in\Z^+$
there are at least $m$ positive integers in the form
$\sum_{s\in I}m_s/n_s$ with $I\se[1,k]$.
Applying this result to the dual $A^*$ of an $m$-system (1.1),
we obtain that there are more than $k-m$ integers in the form
$\sum_{s=1}^kx_s/n_s$ with $x_s\in[0,n_s)$; equivalently, at most
$m-1$ of the numbers in $[1,k]$ cannot be written in the form
$\sum_{s=1}^km_s/n_s=k-\sum_{s=1}^k(n_s-m_s)/n_s$ with $m_s\in[1,n_s]$.
This implies the following result stated in Remark 1.3 of [S03a]:
If (1.1) is an $m$-system, then there are
$m_1,\ldots,m_k\in\Z^+$ such that $\sum_{s=1}^km_s/n_s=m$.

  Our following theorem
 unifies and generalizes many known results.

\proclaim{Theorem 1.1} Let $\Cal A=\{a_s(n_s)\}_{s=0}^k$ be a finite system of
residue classes with $m(\Cal A)>m=\lfloor\sum_{s=1}^km_s/n_s\rfloor$,
where $m_1,\ldots,m_k\in\Z^+$. Then, for any $0\ls \al<1$,
either
$$\sum\Sb I\se[1,k]\\\sum_{s\in I}m_s/n_s=(\al+a)/n_0\endSb
(-1)^{|I|}e^{2\pi i\sum_{s\in I}a_sm_s/n_s}=0\tag1.4$$
for any $a\in\N$, or
$$\bg|\bg\{I\se[1,k]:\,   \sum_{s\in I}\f{m_s}{n_s}=\f{\al+a}{n_0}
\bg\}\bg| \gs\bi m{\lfloor a/n_0\rfloor} \tag1.5$$
for all $a=0,1,2,\ldots$.
\endproclaim

\noindent{\it Example} 1.1. Erd\H os
observed that $\{0(2),0(3),1(4),5(6),7(12)\}$
is a cover of $\Z$ with the moduli
$$n_0=2,\ n_1=3,\ n_2=4,\ n_3=6,\ n_4=12$$
distinct.
As $\lfloor\sum_{s=1}^41/n_s\rfloor=0$,
by Theorem 1.1 in the case $\al=0$ we have
$\sum_{s\in I}1/n_s=1/n_0=1/2$ for some $I\se[1,4]$;
we can actually take $I=\{1,3\}$.
Since $\sum_{s=1}^41/n_s<(5/6+1)/n_0=11/12$,
by Theorem 1.1 in the case $\al=5/6$ the set
$\Cal I=\{I\se[1,4]:\, \sum_{s\in I}1/n_s=5/12\}$
cannot have a single element; in fact,
$\Cal I=\{\{1,4\},\{2,3\}\}$ and
$$(-1)^{|\{1,4\}|}e^{2\pi i(0/n_1+7/n_4)}
+(-1)^{|\{2,3\}|}e^{2\pi i(1/n_2+5/n_3)}
=-e^{\pi i/6}+e^{\pi i/6}=0.$$

\proclaim{Corollary 1.1} If $\Cal A=\{a_s(n_s)\}_{s=0}^k$ is a finite system of
residue classes with $w_{\Cal A}(x)>m=\lfloor\sum_{s=1}^k1/n_s\rfloor$ for all $x\in\Z$, then
$$\bg|\bg\{I\se[1,k]:\, \sum_{s\in I}\f{1}{n_s}=\f{a}{n_0}
\bg\}\bg| \gs\bi m{\lfloor a/n_0\rfloor} \ \ \t{for all}\ a\in\N.\tag1.6$$
In particular, if $(1.1)$ has covering multiplicity
$m(A)=\lfloor\sum_{s=1}^k1/n_s\rfloor$, then
$$\bg|\bg\{I\se[1,k]:\ \sum_{s\in I}\f1{n_s}=n\bg\}\bg|\gs\bi{m(A)}n\ \ \t{for each}\ n\in\N.\tag1.7$$
\endproclaim
\Proof. Observe that the left hand side of (1.4) is nonzero in the case $\al=a=0$.
So (1.6) follows from Theorem 1.1 immediately.
In the case $n_0=1$ this yields the latter result in Corollary 1.1.
\qed
\medskip

\Remark\ 1.1. Let $(1.1)$ be an exact $m$-cover of $\Z$.
Then $\sum_{s=1}^k1/n_s=m$ and $\lfloor\sum_{s\in[1,k]\sm\{t\}}1/n_s\rfloor=m-1$
for any $t=1,\ldots,k$. So Corollary 1.1 implies the following result in [S97]:
For any $t\in[1,k]$ and $a\in\N$, we have
$$\bg|\bg\{I\se[1,k]\sm\{t\}:\ \sum_{s\in I}\f1{n_s}=\f
a{n_t}\bg\}\bg|\gs\bi{m-1}{\lfloor a/n_t\rfloor}.$$
As $m(A)=\sum_{s=1}^k1/n_s$, we also have
$|\{I\se[1,k]:\,   \sum_{s\in I}1/n_s=n\}|\gs \bi mn$ for
all $n=0,1,\ldots,m$, which was first established in [S92b]
by means of the Riemann zeta function.
\endrem

 \proclaim{Corollary 1.2} Let $(1.1)$ be an $m$-system with
$m=\lceil\sum_{s=1}^k1/n_s\rceil$, where $\lceil\al\rceil$ denotes
the least integer not smaller than a real number $\al$. Then
$$\bg|\bg\{\langle m_1,\ldots,m_k\rangle\in\Z^k:\ m_s\in[1,n_s],
\ \sum_{s=1}^k\f{m_s}{n_s}=n\bg\}\bg|\gs\bi{k-m}{n-m}\tag1.8$$ for
every $n=m,\ldots,k$.
\endproclaim
\Proof. Let $n\in[m,k]$. Clearly the left hand side of (1.8) coincides with
$$L:=\bg|\bg\{\langle x_1,\ldots,x_k\rangle:\ x_s\in[0,n_s-1],
\ \sum_{s=1}^k\f{x_s}{n_s}=\sum_{s=1}^k\f{n_s}{n_s}-n=k-n\bg\}\bg|.$$
Since $\sum_{s=1}^k1/n_s>m-1$, $w_A(x)=m$ for some $x\in\Z$.
As the dual $A^*$ of (1.1) has covering multiplicity
$m(A^*)=k-m$, applying Corollary 1.1 to $A^*$ we
find that $L\gs\bi{k-m}{k-n}=\bi{k-m}{n-m}$. This concludes the proof.\qed

\medskip
\Remark\ 1.2. When (1.1) is an exact $m$-cover of $\Z$, it was proved in [S97]
(by a different approach)
that for each $n\in\N$ the equation $\sum_{s=1}^kx_s/n_s=n$ with $x_s\in[0,n_s)$
has at least $\bi{k-m}n$ solutions.

\proclaim{Corollary 1.3}  Let $\Cal A=\{a_s(n_s)\}_{s=0}^k$ be a finite system of
residue classes with $m(\Cal A)>m=\lfloor\sum_{s=1}^km_s/n_s\rfloor$,
where $m_1,\ldots,m_k\in\Z^+$.
Suppose that $J\se[1,k]$ and $\sum_{s\in I}m_s/n_s=\sum_{s\in
J}m_s/n_s$ for no $I\se[1,k]$ with $I\not=J$. Then
$$\bg\{n_0\sum_{s\in J}\f{m_s}{n_s}\bg\}
+\bg\{n_0\sum_{s\in \bar J}\f{m_s}{n_s}\bg\}<1,\tag1.9$$
where $\bar J=[1,k]\sm J$. Also,
$$\sum_{s\in J}\f{m_s}{n_s}\gs m
\ \ \t{or}\ \ \sum_{s\in\bar J}\f{m_s}{n_s}\gs m\tag1.10$$
\endproclaim
\Proof. Let $v=\sum_{s\in J}m_s/n_s$, $\al=\{n_0v\}$ and $b=\lfloor n_0v\rfloor$.
Then $(\al+b)/n_0=v$ and
$$\sum\Sb I\se[1,k]\\\sum_{s\in
I}m_s/n_s=v\endSb(-1)^{|I|}e^{2\pi i\sum_{s\in I}a_sm_s/n_s}
=(-1)^{|J|}e^{2\pi i\sum_{s\in J}a_sm_s/n_s}\not=0.$$
By Theorem 1.1, (1.5) holds for any $a\in\N$.
 Applying (1.5) with $a=mn_0+n_0-1$ we find that
 $\sum_{s\in I}m_s/n_s=(\al+mn_0+n_0-1)/n_0$ for some $I\se[1,k]$,
 therefore $\sum_{s=1}^km_s/n_s\gs m+(\al+n_0-1)/n_0$.
 As $\lfloor \sum_{s=1}^km_s/n_s\rfloor=m$, we must have
$$\bg\{\sum_{s=1}^k\f{m_s}{n_s}\bg\}\gs\f{\al+n_0-1}{n_0},
\ \t{i.e.,}\ n_0-1+\al\ls n_0\bg\{\sum_{s=1}^k\f{m_s}{n_s}\bg\}<n_0.$$
Therefore
$\al\ls\{n_0\{\sum_{s=1}^km_s/n_s\}\}=\{n_0\sum_{s=1}^km_s/n_s\}$,
which is equivalent to (1.9).

(1.5) in the case $a=b$ gives that
 $\bi m{\lfloor b/n_0\rfloor}\ls 1$, thus $\lfloor
v\rfloor\in\{0,m\}$.
As $n_0\{v\}-\al=\lfloor n_0\{v\}\rfloor\ls n_0-1$,
$\{v\}\ls(\al+n_0-1)/n_0\ls\{\sum_{s=1}^km_s/n_s\}$.
If $\lfloor v\rfloor=0$, then $m+v\ls m+\{\sum_{s=1}^km_s/n_s\}=
\sum_{s=1}^km_s/n_s$ and hence
$\sum_{s\in\bar J}m_s/n_s\gs m$.
Therefore (1.10) is valid. We are done. \qed
\medskip

\Remark\ 1.3. Let (1.1) be an exact $m$-cover of $\Z$.
Theorem 4(ii) in [S95] asserts that if $\em\not=J\subset[1,k]$
then $\sum_{s\in I}1/n_s=\sum_{s\in J}1/n_s$ for some $I\se[1,k]$ with
$I\not=J$. This follows from Corollary 1.3, for,
$\Cal A=\{a_s(n_s)\}_{s=0}^k$ (where $a_0=0$ and $n_0=1$) is an
$m+1$-cover of $\Z$ with $\sum_{s\in J\cup\bar J}1/n_s=\sum_{s=1}^k1/n_s=m$.
\endrem

In the 1960s Erd\H os made the following conjecture: For any system (1.1) with $1<n_1<\cdots<n_k$,
if it is a cover of $\Z$ then $\sum_{s=1}^k1/n_s>1$, in other words it cannot be a disjoint cover of $\Z$.
This was later confirmed by H. Davenport, L. Mirsky, D. Newman and R. Rad\'o who proved that
if (1.1) is a disjoint cover of $\Z$ with $1<n_1\ls\cdots\ls n_{k-1}\ls n_k$ then $n_{k-1}=n_k$.

\proclaim{Corollary 1.4} Let $(1.1)$ be an $m$-cover of $\Z$ with
$$n_1\ls\cs\ls n_{k-l}<n_{k-l+1}=\cs=n_k\ \ (0<l<k).\tag1.11$$
Then, for any $r\in[0,l]$ with $r<n_k/n_{k-l}$,
either $\sum_{s=1}^{k-r}1/n_s\gs m$ or
$$\bi lr\in D(n_k)=\bg\{\sum_{p\mid n_k}px_p:\ x_p\in\N
\ \t{for any prime divisor}\ p\ \t{of}\ n_k\bg\}.$$
\endproclaim
\Proof. Set $\Cal A=\{a_s(n_s)\}_{s=0}^k$ where $a_0=0$ and $n_0=1$.
Suppose that $\sum_{s=1}^{k-r}1/n_s<m$. Then $\sum_{s=1}^k1/n_s<m+r/n_k<m+1\ls m(\Cal A)$.
Since $|\{I\se[1,k]:\, \sum_{s\in I}1/n_s=m+r/n_k\}|=0<\bi mm$, by Theorem 1.1 we must have
$$\sum\Sb I\se[1,k]\\\sum_{s\in I}1/n_s=r/n_k\endSb
(-1)^{|I|}e^{2\pi i\sum_{s\in I}a_sm_s/n_s}=0.$$
Observe that $r/n_k<1/n_{k-l}=\min\{1/n_s:\, 1\ls s\ls k-l\}$. Therefore
$$0=\sum\Sb I\se(k-l,k]\\\sum_{s\in I}1/n_s=r/n_k\endSb
(-1)^{|I|}e^{2\pi i\sum_{s\in I}a_s/n_s}
=(-1)^r\Sigma_r,$$
where
$$\Sigma_r=\sum\Sb I\se(k-l,k]\\|I|=r\endSb
e^{2\pi i\sum_{s\in I}a_s/n_k}.$$
By Lemma 3.1 of [S03a],  $\Sigma_r=0$ implies
that $$\bi lr=|\{I\se(k-l,k]:\,  |I|=r\}|\in D(n_k).$$
This concludes the proof. \qed
\medskip

\Remark\ 1.4. Let (1.1) be an $m$-cover of $\Z$ with (1.11).
By Corollary 1.4 in the case $r=l$, either $l\gs n_k/n_{k-l}>1$
or $\sum_{s=1}^{k-l}1/n_s\gs m$; this is one of the main results in [S96].
Corollary 1.4 in the case $r=1$ yields that either $\sum_{s=1}^{k-1}1/n_s\gs m$
or $l\in D(n_k)$; this implies the extended Newman-Zn\'am result (cf. [N71]) which asserts that
if (1.1) is an exact $m$-cover of $\Z$ (and hence $\sum_{s=1}^{k-1}1/n_s<\sum_{s=1}^k1/n_s=m$)
then $l$ is not smaller than the least prime divisor of $n_k$.
\endrem

Let $(1.1)$ is an $m$-system with (1.11), and let $r\in\N$ and $r<n_k/n_{k-l}$.
With the help of the dual system of (1.1), we can also show that either
 $\sum_{s=1}^k1/n_s\ls m-r/n_k$ or
 $$\bi{l+r-1}r=|\{\langle x_{k-l+1},\ldots,x_k\rangle\in\N^l:\ x_{k-l+1}+\cdots+x_k=r\}|\in D(n_k).$$

 If (1.1) is disjoint with $1<n_1<\cs<n_k$,
 then $\sum_{s=1}^k1/n_s<1$ since (1.1) is not a disjoint cover of $\Z$;
Erd\H os [E62] showed further that $\sum_{s=1}^k1/n_s\ls 1-1/2^k$.
Now we give a generalization of this result.

 \proclaim{Theorem 1.2} Let $(1.1)$ be an $m$-system with $k\gs m$,
$\sum_{s=1}^k1/n_s\not=m$ and $n_1\ls\cs\ls n_k$. Then we have
$$\sum_{s=1}^k\f1{n_s}\ls m-\f1{2^{k-m+1}},\tag1.12$$
and equality holds if and only if $n_s=2^{\max\{s-m+1,0\}}$ for
all $s=1,\ldots,k$.
\endproclaim
\Remark\ 1.5. Let $k\gs m\gs1$ be integers. Then $m-1$ copies of
$0(1)$, together with the following $k-m+1$ residue classes
$$1(2),\ 2(2^2),\ \ldots,\ 2^{k-m}(2^{k-m+1}),$$
form an $m$-system with the moduli $2^{\max\{s-m+1,0\}}\
(s=1,\ldots,k)$.
\endrem

  We will prove Theorems 1.1 and 1.2 in the next section.
Section 3 deals with two characterizations of $m$-systems one
of which is as follows.

\proclaim{Theorem 1.3}  $(1.1)$ is an $m$-system if and only if
for any $n\in[m,k)$ we have
$$S(n,\al)=\cases(-1)^k\ &\t{if}\ \al=0,\\0\
&\t{if}\ 0<\al<1,\endcases\tag1.13$$
where $S(n,\al)$ represents the sum
$$\sum\Sb m_1,\ldots,m_k\in\Z^+\\\{\sum_{s=1}^km_s/n_s\}=\al\endSb
(-1)^{\lfloor\sum_{s=1}^km_s/n_s\rfloor}\bi
n{\lfloor\sum_{s=1}^km_s/n_s\rfloor} e^{2\pi i\sum_{s=1}^ka_sm_s/n_s}.$$
\endproclaim

Theorem 1.3 in the case $m=1$ yields the following result.
\proclaim{Corollary 1.5} If $(1.1)$ is disjoint, then we have
$$\sum\Sb m_1,\ldots,m_k\in\Z^+\\\sum_{s=1}^k m_s/n_s=1\endSb
 e^{2\pi i\sum_{s=1}^ka_sm_s/n_s}=(-1)^{k-1}.
\tag1.14$$
\endproclaim

We conclude this section with the following challenging conjecture
arising from the author's solutions (cf. [S92a])
to the two open problems posed in [HM].
\proclaim{Conjecture 1.1} For any disjoint system $(1.1)$ with $k>1$,
there exists a pair $\{s,t\}\ (1\ls s<t\ls k)$ such that $\gcd(n_s,n_t)\gs k$.
\endproclaim

\Remark\ 1.6. In 1992 the author [S92a] verified the conjecture for $k\ls 4$.

\heading{2. Proofs of Theorems 1.1 and 1.2}\endheading

\proclaim{Lemma 2.1} Let $N\in\Z^+$ be a common multiple of the moduli $n_1,\ldots,n_k$
in $(1.1)$. And let $m,m_1,\ldots,m_k\in\Z^+$. If $(1.1)$ is an $m$-cover
of $\Z$, then $(1-z^N)^m$ divides the polynomial
$\prod_{s=1}^k(1-z^{Nm_s/n_s}e^{2\pi ia_sm_s/n_s})$.
When $m_1,\ldots,m_k$ are relatively prime to $n_1,\ldots,n_k$
respectively, the converse also holds.
\endproclaim
\Proof. For any $r=0,1,\ldots,N-1$, clearly $e^{2\pi i r/N}$
 is a zero of the polynomial $\prod_{s=1}^k(1-z^{Nm_s/n_s}e^{2\pi ia_sm_s/n_s})$
 with multiplicity $M_r=|\{s\in[1,k]:\, n_s\mid m_s(r+a_s)\}|$.
 Observe that $M_r\gs w_A(-r)$. If $m_s$ is relatively prime to
 $n_s$ for each $s\in[1,k]$,
 then $M_r=w_A(-r)$.
 As $(1-z^N)^m=\prod_{r=0}^{N-1}(1-ze^{-2\pi i r/N})^m$,
 the desired result follows from the above.

\medskip
\noindent{\it Proof of Theorem 1.1}. Set $m_0=1$, and let $N_0$
be the least common multiple of $n_0,n_1,\ldots,n_k$. In light of Lemma 2.1, we can write
$P(z)=\prod_{s=0}^k(1-z^{N_0m_s/n_s}e^{2\pi ia_sm_s/n_s})$ in the form
$(1-z^{N_0})^{m+1}Q(z)$ where $Q(z)\in\C[z]$. Clearly
$$\deg Q=\deg P-(m+1)N_0=N_0\(\sum_{s=0}^k\f{m_s}{n_s}-m-1\)<\f {N_0}{n_0}.$$
Also,
$$\aligned&\prod_{s=1}^k\l(1-z^{N_0m_s/n_s}e^{2\pi ia_sm_s/n_s}\r)
\\=&\sum_{n=0}^m(-1)^n\bi mnz^{nN_0}\sum_{r=0}^{n_0-1}z^{rN_0/n_0}e^{2\pi i ra_0/n_0}Q(z)
\endaligned\tag2.1$$
since
$$\f{1-z^{N_0}}{1-z^{N_0/n_0}e^{2\pi ia_0/n_0}}
=\sum_{r=0}^{n_0-1}z^{rN_0/n_0}e^{2\pi i ra_0/n_0}.$$

Let $a\in\N$ and
$$C_a=(-1)^{\lfloor a/n_0\rfloor}
\sum\Sb I\se[1,k]\\\sum_{s\in I}m_s/n_s=(\al+a)/n_0\endSb
(-1)^{|I|}e^{2\pi i\sum_{s\in I}(a_s-a_0)m_s/n_s}.$$ By comparing the
coefficients of $z^{N_0(\al+a)/n_0}$ on both sides of (2.1) we
obtain that
$$\align&\sum\Sb I\se[1,k]\\\sum_{s\in I}m_s/n_s=(\al+a)/n_0\endSb
(-1)^{|I|}e^{2\pi i\sum_{s\in I}a_sm_s/n_s}
\\=&(-1)^{\lfloor a/n_0\rfloor}\bi m{\lfloor a/n_0\rfloor}e^{2\pi ia_0\{a/n_0\}}
[z^{\al N_0/n_0}]Q(z),
\endalign$$
where $[z^{\al N_0/n_0}]Q(z)$ denotes the coefficient of
$z^{\al N_0/n_0}$ in $Q(z)$.
Therefore
$$C_a=e^{-2\pi i\al a_0/n_0}\bi m{\lfloor a/n_0\rfloor}[z^{\al N_0/n_0}]Q(z)
=\bi m{\lfloor a/n_0\rfloor}C_0.\tag2.2$$

For an algebraic integer
$\omega$ in the field $K=\Q(e^{2\pi i/N_0})$,
the norm $N(\omega)=\prod_{1\ls r\ls N_0,\ \gcd(r,N_0)=1}\sigma_r(\omega)$
(with respect to the field
extension $K/\Q$) is a rational integer, where $\sigma_r$
is the automorphism of $K$ (in the Galois group $\t{Gal}(K/\Q)$)
induced by $\sigma_r(e^{2\pi i/N_0})=e^{2\pi ir/N_0}$.
(See, e.g., [K97].)
As $N((-1)^{\lfloor a/n_0\rfloor}C_a)$ equals
$$\prod\Sb 1\ls r\ls N_0\\\gcd(r,N_0)=1\endSb
\ \sum\Sb I\se[1,k]\\\sum_{s\in I}m_s/n_s=(\al+a)/n_0\endSb
(-1)^{|I|}e^{2\pi i r\sum_{s\in I}(a_s-a_0)m_s/n_s},$$
we have
$$\align |N(C_a)|=&\prod\Sb 1\ls r\ls N_0\\\gcd(r,N_0)=1\endSb
\bg|\sum\Sb I\se[1,k]\\\sum_{s\in I}m_s/n_s=(\al+a)/n_0\endSb
(-1)^{|I|}e^{2\pi i r\sum_{s\in I}(a_s-a_0)m_s/n_s}\bg|
\\\ls&\bg|\bg\{I\se[1,k]:\,   \sum_{s\in I}\f{m_s}{n_s}
=\f{\al+a}{n_0}\bg\}\bg|^{\varphi(N_0)},
\endalign$$
where $\varphi$ is Euler's totient function.
Also,
$$|N(C_a)|=\bg|N\l(\bi m{\lfloor
a/n_0\rfloor}\r)\bg|\times|N(C_0)|
=\bi m{\lfloor a/n_0\rfloor}^{\varphi(N_0)}|N(C_0)|.$$

Suppose that $C_b\not=0$ for some $b\in\N$. Then $N(C_b)\not=0$
and hence $N(C_0)\in\Z$ is nonzero. For any $a\in\N$, we have
$$\bg|\bg\{I\se[1,k]:\, \sum_{s\in I}\f{m_s}{n_s}
=\f{\al+a}{n_0}\bg\}\bg|^{\varphi(N_0)}
\gs |N(C_a)|\gs\bi m{\lfloor a/n_0\rfloor}^{\varphi(N_0)}$$
and hence (1.5) holds. This concludes the proof. \qed

\medskip
\noindent{\it Proof of Theorem 1.2}. We use induction on $k$.

In the case $k=m$,
we have $n_k>1$ and hence
$$\sum_{s=1}^k\f1{n_s}\ls k-1+\f1{n_k}\ls m-\f12=m-\f1{2^{k-m+1}},$$
also $\sum_{s=1}^k1/n_s=m-1/2$ if and only if $n_1=\cs=n_{k-1}=1$
and $n_k=2$.

Now let $k>m$. Clearly $\sum_{s=1}^{k-1}1/n_s<\sum_{s=1}^k1/n_s<m$. Assume that
$$\sum_{s=1}^{k-1}\f1{n_s}\ls m-\f1{2^{(k-1)-m+1}}=m-\f1{2^{k-m}}$$
and that equality holds if and only if $n_s=2^{\max\{s-m+1,0\}}$
for all $s\in[1,k-1]$.
When $n_k>2^{k-m+1}$, we have
$$\sum_{s=1}^k\f1{n_s}=\sum_{s=1}^{k-1}\f1{n_s}+\f1{n_k}
<\l(m-\f1{2^{k-m}}\r)+\f1{2^{k-m+1}}=m-\f1{2^{k-m+1}}.$$

If $\sum_{s=1}^k1/n_s>m-1/n_k$, then
$\lceil\sum_{s=1}^k1/n_s\rceil=m$,
thus $\sum_{s=1}^km_s/n_s=m$ for some $m_1,\ldots,m_k\in\Z^+$
(by Corollary 1.2) and hence
$$m-\sum_{s=1}^k\f1{n_s}\gs\min\l\{\f1{n_s}:1\ls s\ls k\r\}=\f1{n_k}.$$
This shows that we do have $\sum_{s=1}^k1/n_s\ls m-1/n_k$.
Providing $n_k\ls2^{k-m+1}$, (1.12) holds, and also
$$\align\sum_{s=1}^k\f1{n_s}=m-\f1{2^{k-m+1}}&\iff n_k=2^{k-m+1}
\ \t{and}\ \sum_{s=1}^{k-1}\f1{n_s}=m-\f1{2^{k-m}}
\\&\iff n_s=2^{\max\{s-m+1,0\}}\ \t{for}\ s=1,\ldots,k-1,k.
\endalign$$
This concludes the induction step and we are done. \qed

\heading{3. Characterizations of $m$-systems}\endheading

\noindent{\it Proof of Theorem 1.3}.
Like Lemma 2.1, (1.1) is an
$m$-system if and only if
$f(z)=(1-z^N)^m/\prod_{s=1}^k(1-z^{N/n_s}e^{2\pi ia_s/n_s})$ is a
polynomial, where $N$ is the least common
multiple of $n_1,\ldots,n_k$.

Set $c=m-\sum_{s=1}^k1/n_s$.  If $f(z)$ is a polynomial, then
$\deg f=cN$ and $[z^{cN}]f(z)=(-1)^{k-m}e^{-2\pi i\sum_{s=1}^ka_s/n_s}$.

  For $|z|<1$ we have
$$f(z)=\sum_{n=0}^m\bi
mn(-1)^nz^{nN}\prod_{s=1}^k\sum_{x_s=0}^{\infty}e^{2\pi i
a_sx_s/n_s}z^{Nx_s/n_s}.$$ Let $\al\gs0$. Then
$$\align[z^{(c+\al)N}]f(z)=&\sum_{n=0}^m(-1)^n\bi mn\sum\Sb
x_1,\ldots,x_k\in\N\\\sum_{s=1}^kx_s/n_s=c+\al-n\endSb e^{2\pi
i\sum_{s=1}^ka_sx_s/n_s}
\\=&\sum_{n=0}^m(-1)^n\bi mn\sum\Sb
m_1,\ldots,m_k\in\Z^+\\\sum_{s=1}^km_s/n_s=\al+m-n\endSb e^{2\pi
i\sum_{s=1}^ka_s(m_s-1)/n_s}
\\=&(-1)^me^{-2\pi i\sum_{s=1}^ka_s/n_s}S(m,\al),\endalign$$
where $S(n,\al)\ (n\in\N)$ represents the sum
$$\sum\Sb
m_1,\ldots,m_k\in\Z^+\\\sum_{s=1}^km_s/n_s-\al\in\N\endSb
(-1)^{\sum_{s=1}^km_s/n_s-\al}\bi n{\sum_{s=1}^km_s/n_s-\al}
e^{2\pi i\sum_{s=1}^ka_sm_s/n_s}$$ which agrees with its
definition in the case $0\ls\al<1$ given in Theorem 1.3.

 (i) Suppose that (1.1) is an $m$-system. Then $f(z)$ is a
polynomial of degree $cN$ and hence
$$S(m,\al)=(-1)^me^{2\pi i\sum_{s=1}^ka_s/n_s}[z^{(c+\al)N}]f(z)
=\cases(-1)^k\ &\t{if}\ \al=0,\\0&\t{if}\ \al>0.\endcases$$ For
any integer $n\gs m$, (1.1) is also an $n$-system and so we have (1.13).

 (ii) Now assume that $(1.13)$ holds for all $n\in[m,k)$.
For any $n\gs k$ we also have (1.13) by (i) because
(1.1) is a $k$-system.

If $0<\al<1$ then $S(n,\al)=0$ for any integer $n\gs m$. Fix
$\al>0$. If $S(n,\al)=0$ for all integers $n\gs m$, then for any
integer $n\gs m$ we have
$$S(n,\al+1)=S(n,\al)-S(n+1,\al)=0$$
because $\bi n{j-1}=\bi{n+1}j-\bi nj$ for $j=1,2,\ldots$.
Thus, by induction,
 $S(n,\al)=0$ for all $\al>0$ and $n=m,m+1,\ldots$.
It follows that $[z^{(c+\al)N}]f(z)=0$ for any $\al>0$. So $f(z)$
is a polynomial and (1.1) is an $m$-system.

The proof of Theorem 1.3 is now complete. \qed

    The following characterization of $m$-covers
 plays important roles in [S95, S96].

\proclaim{Lemma 3.1 {\rm ([S95])}} Let $m,m_1,\ldots,m_k\in\Z^+$. If
$(1.1)$ forms  an $m$-cover of $\Z$, then
$$\sum\Sb I\se[1,k]\\\{\sum_{s\in I}m_s/n_s\}=\theta\endSb
(-1)^{|I|}\bi{\lfloor\sum_{s\in I}m_s/n_s\rfloor}{n}e^{2\pi
i\sum_{s\in I}a_sm_s/n_s} =0\tag3.1$$ for all $0\ls\theta<1$ and
$n=0,1,\ldots,m-1$. We also have the converse if $m_1,\ldots,m_k$ are
relatively prime to $n_1,\ldots,n_k$ respectively.
\endproclaim

 We can provide a new proof of Lemma 3.1
 in a way similar to the proof of Theorem 1.3.

\proclaim{Lemma 3.2} Let $n\in\Z^+$ and
$l\in[0,n-1]$. Then
$$\sum\Sb J\se[1,n)\\|J|=l\endSb e^{2\pi i\sum_{j\in
J}j/n}=(-1)^l.\tag3.2$$
\endproclaim
\Proof. Clearly we have the identity
$$\prod_{0<j<n}\l(1-ze^{2\pi ij/n}\r)=\f{1-z^n}{1-z}=1+z+\cs+z^{n-1}.$$
Comparing the coefficients of $z^l$ we then obtain (3.2). \qed

Using Lemmas 3.1 and 3.2 we can deduce another characterization of
$m$-systems.

\proclaim{Theorem 3.1} $(1.1)$ is an $m$-system
 if and only if we have
$$\sum\Sb x_s\in [0,n_s)\ \t{for}\ s\in[1,k]
\\\{\sum_{s=1}^kx_s/n_s\}=\theta\endSb
\bi{\lfloor\sum_{s=1}^kx_s/n_s\rfloor}ne^{2\pi i\sum_{s=1}^k
a_sx_s/n_s}=0\tag3.3$$
for all $0\ls\theta<1$ and $n\in[0,k-m)$.
\endproclaim
\Proof. The case $k\ls m$ is trivial, so we just let $k>m$.
Recall that (1.1) is an $m$-system if and only if its dual $A^*$
is a $(k-m)$-cover of $\Z$.

By Lemma 3.1 in the case $m_1=\cs=m_k=1$, $A^*$ forms an
$(k-m)$-cover of $\Z$ if and only if for any $0\ls\theta<1$ and
$n\in[0,k-m)$ the sum
$$\sum\Sb x_s\in [0,n_s)\ \t{for}\ s\in[1,k]
\\\{\sum_{s=1}^kx_s/n_s\}=\theta\endSb(-1)^{\sum_{s=1}^kx_s}
\bi{\lfloor\sum_{s=1}^kx_s/n_s\rfloor}ne^{2\pi i\sum_{s=1}^k
a_sx_s/n_s}\prod_{s=1}^kf_s(x_s)$$ vanishes, where
$$f_s(x_s)=\sum\Sb J\se[1,n_s)\\|J|=x_s\endSb e^{2\pi i\sum_{j\in
J}j/n_s}=(-1)^{x_s}$$ by Lemma 3.2.
This concludes the proof. \qed

The following consequence extends Corollary 1.5.
\proclaim{Corollary 3.1} Let $(1.1)$ be an $m$-system. Then we
have
$$\sum\Sb  m_s\in[1,n_s]\ \t{for}\ s\in[1,k]\\m-\sum_{s=1}^k m_s/n_s\in\N\endSb
\bi{k-1-\sum_{s=1}^km_s/n_s}{m-\sum_{s=1}^km_s/n_s} e^{2\pi
i\sum_{s=1}^ka_sm_s/n_s}=(-1)^{k-m}.$$
\endproclaim
\Proof.  If $k\ls m$, then the left hand side of the last equality
coincides with
$$\bi{k-1-\sum_{s=1}^kn_s/n_s}{m-\sum_{s=1}^kn_s/n_s}
e^{2\pi i\sum_{s=1}^ka_sn_s/n_s}=\bi{-1}{m-k}=(-1)^{m-k}.$$

Now let $k>m$. As $\{-a_s(n_s)\}^k_{s=1}$ is an $m$-system, by
Theorem 4.1 and the identity
$$(-1)^{k-m-1}\bi{x-1}{k-m-1}=\sum_{n=0}^{k-m-1}(-1)^n\bi xn$$
(cf. [GKP, (5.16)]) we have
$$\align0=&\sum\Sb x_s\in[0,n_s)\ \t{for}\ s\in[1,k]
\\\{\sum_{s=1}^kx_s/n_s\}=0\endSb\bi{\lfloor\sum_{s=1}^kx_s/n_s\rfloor-1}{k-m-1}
e^{2\pi i\sum_{s=1}^k(-a_s)x_s/n_s}
\\=&\sum\Sb m_s\in[1,n_s]\ \t{for}\ s\in[1,k]\\\sum_{s=1}^k(n_s-m_s)/n_s\in\N\endSb
\bi{\sum_{s=1}^k(n_s-m_s)/n_s-1}{k-m-1}e^{-2\pi
i\sum_{s=1}^ka_s(n_s-m_s)/n_s}
\\=&\sum\Sb m_s\in[1,n_s]\ \t{for}\ s\in[1,k]\\\sum_{s=1}^km_s/n_s\in[0,k-1]\endSb
\bi{k-1-\sum_{s=1}^km_s/n_s}{k-1-m}e^{2\pi
i\sum_{s=1}^ka_sm_s/n_s}
\\&\qquad+\bi{k-1-\sum_{s=1}^kn_s/n_s}{k-1-m}e^{2\pi i\sum_{s=1}^ka_sn_s/n_s}.
\endalign$$
So the desired equality follows. \qed


\widestnumber\key{GKP}
\Refs

\ref\key C71\by R. Crocker\paper On a sum of a prime and two
powers of two \jour Pacific J.
Math.\vol36\yr1971\pages103--107\endref

\ref\key E50\by P. Erd\H os\paper On integers of the form $2^k+p$
and some related problems\jour Summa Brasil.
Math.\vol2\yr1950\pages113--123\endref

\ref\key E62\by P. Erd\H os\paper Remarks on number theory {\rm
IV:} Extremal problems in number theory {\rm I}\jour Mat.
Lapok\vol 13\yr1962\pages228-255\endref

\ref\key F02\by M. Filaseta\paper
Coverings of the integers associated with an irreducibility theorem
of A. Schinzel \jour in: Number Theory for the Millennium
(Urbana, IL, 2000), Vol. II, pp. 1-24, A K Peters, Natick, MA,
2002\endref

\ref\key GKP\by R. L. Graham, D. E. Knuth and O. Patashnik
\book Concrete Mathematics: A Foundation for Computer
Science {\rm (2nd edition)} \publ Addison-Wesley, Amsterdam \yr 1994\endref

\ref\key HM\by A. P. Huhn and L. Megyesi\paper On disjoint residue classes
\jour Discrete Math.\vol 41\yr 1982\pages 327--330\endref

\ref\key K97\by H. Koch\book Algebraic Number Theory
\publ Springer, Berlin, 1997, Chapter 1\endref

\ref\key N71\by M. Newman\paper Roots of unity and covering
sets\jour Math. Ann.\vol 191\yr 1971\pages 279--282\endref

\ref\key P76\by \v S. Porubsk\'y\paper On $m$ times covering
systems of congruences\jour Acta Arith.\vol
29\yr1976\pages159--169\endref

\ref\key PS\by\v S. Porubsk\'y and J. Sch\"onheim \paper Covering
systems of Paul Erd\"os: past, present and future \jour in: Paul
Erd\"os and his Mathematics. I (edited by G. Hal\'asz, L.
Lov\'asz, M. Simonvits, V. T. S\'os), Bolyai Soc. Math. Studies
11, Budapest, 2002, pp. 581--627\endref

\ref\key S92a\by Z. W. Sun\paper Solutions to two problems of Huhn and Megyesi
\jour Chinese Ann. Math. Ser. A\vol 13\yr 1992\pages 722--727\endref

\ref\key S92b\by Z. W. Sun\paper On exactly $m$ times covers \jour
Israel J. Math. \vol 77\yr 1992\pages 345--348\endref

\ref\key S95\by Z. W. Sun\paper Covering the integers by
arithmetic sequences
 \jour Acta Arith.\vol 72\yr1995\pages109--129\endref

\ref\key S96\by Z. W. Sun\paper Covering the integers by
arithmetic sequences {\rm II}
 \jour Trans. Amer. Math. Soc.\vol348\yr1996\pages4279--4320\endref

\ref\key S97\by Z. W. Sun\paper Exact $m$-covers and the linear
form $\sum^k_{s=1} x_s/n_s$ \jour Acta Arith.\vol81\yr1997\pages
175--198\endref

\ref\key S00\by Z. W. Sun\paper On integers not of the form $\pm
p^a\pm q^b$ \jour Proc. Amer. Math. Soc.\vol 128\yr 2000\pages
997--1002\endref

\ref\key S01\by Z. W. Sun\paper Algebraic approaches to periodic
arithmetical maps\jour J. Algebra\vol 240\yr
2001\pages723--743\endref

\ref\key S03a\by Z. W. Sun\paper On the function $w(x)=|\{1\le
s\le k:\, x\eq a_s\ (\mo\ n_s)\}|$ \jour Combinatorica\vol
23\yr 2003\pages 681--691\endref

\ref\key S03b\by Z. W. Sun\paper Unification of zero-sum problems,
subset sums and covers of $\Z$ \jour Electron. Res. Annnounc.
Amer. Math. Soc. \vol 9\yr 2003\pages 51--60\endref

\ref\key S04\by Z. W. Sun\paper Arithmetic properties of periodic maps
\jour Math. Res. Lett.\vol 11\yr 2004\pages to appear\endref

\ref\key SY\by Z. W. Sun and S. M. Yang\paper A note on integers
of the form $2^n+cp$ \jour Proc. Edinburgh. Math. Soc.\vol 45\yr
2002 \pages 155--160\endref

\endRefs
\enddocument